\documentclass[12pt,a4paper]{article}
\usepackage{amsmath}
\usepackage{amsthm}
\usepackage{amssymb}
\usepackage{amsbsy}
\usepackage{amsfonts}
\usepackage{amstext}
\usepackage{amscd}
\usepackage[dvips]{graphicx,psfrag}
\numberwithin{equation}{section}
\theoremstyle{plain}
\newtheorem{thm}{Theorem}[section]
\newtheorem{prop}[thm]{Proposition}

\newtheorem{lem}[thm]{Lemma}
\theoremstyle{definition}
\newtheorem{exa}[thm]{Example}

\newtheorem{rem}[thm]{Remark}
\newtheorem{defi}[thm]{Definition}

\newcommand{\pa}{\mathcal{P}}

\newcommand{\real}{\mathbb{R}}
\newcommand{\rat}{\mathbb{Q}}
\newcommand{\tor}{\mathbb{T}}
\newcommand{\disc}{\mathbb{D}}
\DeclareMathOperator*{\com+}{\mathbb{C}_+}
\newcommand{\comp}{\mathbb{C}}
\newcommand{\nat}{\mathbb{N}}
\newcommand{\im}{\text{\normalfont Im}\,}
\newcommand{\re}{\text{\normalfont Re}\,}

\def\mrhd{\kern0.2em\rule{0.035em}{0.52em}\kern-.35em\gtrdot \kern-.2em}
\def\submrhd{\mathop{\kern0.1em\lower0.1ex\hbox{\rule{0.03em}{0.42em}}\kern-.1em\gtrdot \kern0.02em}}
\def\dismrhd{\mathop{{\kern0.1em\lower0.07ex\hbox{\rule{0.035em}{0.57em}}\kern-.1em\gtrdot \kern0.02em}}}

\newcommand{\autimes}{\!\!\begin{array}{c} {\scriptstyle\times} \\[-12pt]\cup\end{array}\!\!}

\newcommand{\utimes}{\kern0.05em\buildrel{\times}\over{\rule{0em}{0.004em}} \kern-0.9em\cup \kern0.2em}
\newcommand{\putimes}{\mathop{\kern0.05em\buildrel{\times}\over{\rule{0em}{0.004em}} \kern-0.9em\cup \kern0.2em}}
\newcommand{\hutimes}{\mathop{\kern0.02em\buildrel{\times}\over{\rule{0em}{0.004em}} \kern-0.48em\cup \kern-0.1em}}
\begin{document}
\title{Conditionally monotone independence II: Multiplicative convolutions and infinite divisibility}
\author{Takahiro Hasebe\footnote{Supported by Grant-in-Aid for JSPS Research Fellows.} \\ Graduate School of Science, Kyoto University, \\ Kyoto 606-8502, Japan \\ E-mail: hsb@kurims.kyoto-u.ac.jp \\ Tel: +81-75-753-7245 \\Fax: +81-75-753-7272 }
\date{}

\maketitle

\begin{abstract}
We study the multiplicative convolution for c-monotone independence. This convolution unifies the monotone, Boolean and orthogonal multiplicative convolutions. We characterize convolution semigroups for the c-monotone multiplicative convolution on the unit circle. 
We also prove that an infinitely divisible distribution can always be embedded in a convolution semigroup. We furthermore discuss the (non)-uniqueness of such embeddings including the monotone case. Finally connections to the multiplicative Boolean 
convolution are discussed.  
\end{abstract}



\section{Introduction}
In non-commutative probability theory, many kinds of independence are known. 
Among them, tensor, free, Boolean and monotone independences \cite{Mur3,S-W,V1} are important since they satisfy natural properties \cite{Mur5}.  
Free, Boolean and monotone independences can be unified in terms of conditionally free (c-free, for simplicity) independence \cite{BLS,Fra}; free cumulants \cite{V1} and Boolean cumulants \cite{S-W} can also be unified by c-free cumulants \cite{BLS}. Only in the monotone case, however, monotone cumulants \cite{H-S} cannot be unified by c-free cumulants. To overcome this difficulty, conditionally monotone (c-monotone, for simplicity) independence has been introduced in \cite{Has3}, and as a result, orthogonal independence \cite{Len1} turned out to be included in c-free independence and c-monotone independence.  


In this paper we study the multiplicative convolution associated to c-monotone independence. While c-monotone cumulants cannot be unified by c-free cumulants, the complex analytic characterization of the additive c-monotone convolution follows from the additive c-free convolution. This situation is the same for multiplicative convolutions: we show a complex analytic characterization of the multiplicative c-monotone convolution by using a result of the c-free case. Then we characterize infinitely divisible distributions. These results can be seen as a generalization of results of papers \cite{Ber1, Fra2}. 

Multiplicative convolutions sometimes cause problems which do not appear in additive convolutions. For instance, while the monotone and orthogonal convolutions preserve the probability measures on $[0,\infty)$, the Boolean convolution does not \cite{Ber1,Ber2,Fra2,Len2}. Another instance is the fact that a probability measure on the unit circle is not always infinitely divisible with respect to the Boolean convolution \cite{Fra2}. This makes it difficult to define a multiplicative analogue of $t$-transformation, an additive version of which was first introduced in \cite{BW1} to deform the additive free convolution. 
The latter problem will be understood more in this paper.  

Let us explain the main contents of each section. 
Section \ref{sec1} is devoted to relations among kinds of multiplicative convolutions in a unified way in terms of c-free convolutions.  
In Section \ref{sec2} we characterize multiplicative c-monotone convolutions by using generating functions for the c-free convolution.   
In Section \ref{sec3} we prove that there exists a one-to-one correspondence between a c-monotone convolution semigroup and a pair of analytic vector fields. 
In Section \ref{sec4} we prove an embedding of  an infinitely divisible distribution into a convolution semigroup. The proof uses results on monotone convolutions of \cite{Ber1}. 
We discuss uniqueness and non-uniqueness of such embeddings.  
We then discuss a multiplicative version of $t$-transformation coming from multiplicative Boolean convolutions.

Notation and necessary concepts are provided below.  
An algebraic probability space is a pair $(\mathcal{A},\varphi)$ of a unital algebra $\mathcal{A}$ and a linear functional $\varphi$. $X \in \mathcal{A}$ is called a random variable. We always assume that a linear functional preserves the unit. If we consider probability distributions of random variables, then positivity is needed in $(\mathcal{A},\varphi)$, so that we assume that $\mathcal{A}$ is a $\ast$- (or  $C^\ast$-) algebra and $\varphi$ is a state. If two linear functionals $\varphi, \psi$ are provided in the same algebra, we also call a triple $(\mathcal{A},\varphi,\psi)$ an algebraic probability space. 

We define two independences which are in particular important in this article. 
C-monotone independence was introduced in \cite{Has3}. 
\begin{defi}
Let $(\mathcal{A}, \varphi, \psi)$ be an algebraic probability space; let $I$ be a linearly ordered set.   
We consider sublagebras $\{\mathcal{A}_i \}_{i \in I}$, each of which does not contain the unit of $\mathcal{A}$.  
 $\mathcal{A}_i$ are said to be c-monotone independent if 
the following properties are satisfied for all elements $X_k \in \mathcal{A}_{i_k}$ and indices $i_1,\cdots,i_n$, $n \geq 1$: 
\begin{itemize}
\item[(1)] $\varphi (X_1 \cdots X_n ) = \varphi (X_1) \varphi(X_2 \cdots X_n)$  \text{~whenever~} $i_1 > i_2$;    
\item[(2)] $\varphi (X_1 \cdots X_n ) = \varphi(X_1\cdots X_{n-1})\varphi(X_n)$ \text{~whenever~} $i_n > i_{n-1}$;   
\item[(3)]  $\varphi (X_1 \cdots X_n)  = (\varphi (X_j) - \psi (X_j)) \varphi (X_1 \cdots X_{j-1})  \varphi (X_{j+1} \cdots X_n) + \psi (X_j) \varphi (X_1\cdots X_{j-1} X_{j+1} \cdots X_n)$ \text{~whenever~} $j$ satisfies $i_{j-1} < i_j > i_{j+1}$ and $2 \leq j \leq n-1$; 
\item[(4)] $\mathcal{A}_i$ are monotone independent with respect to $\psi$.  
\end{itemize} 
\end{defi}
Monotone independence was defined in \cite{Mur3} (see also \cite{Mur2}). We however note that the above properties (1)-(3) become monotone independence with respect to $\psi$ in the special case $\varphi = \psi$.  

C-monotone independence was defined for subalgebras; however independence for random variables $\{X_i \}$ can also be defined if we consider the subalgebra $\mathcal{A}_i$ generated by $X_i$ without the unit of $\mathcal{A}$. 

C-free independence was introduced in \cite{BS} and further studied in \cite{BLS}. For later use, we also 
define a c-free product of algebraic probability spaces.  
\begin{defi} (1) Let $I$ be an index set and let $(\mathcal{A}_i, \varphi_i, \psi_i)$ be algebraic probability spaces. 
The c-free product $(\mathcal{A}, \varphi, \psi) = \ast_{i \in I}(\mathcal{A}_i, \varphi_i, \psi_i)$ is defined as follows. 
$\mathcal{A} := \ast_{i \in I}\mathcal{A}_i$ is the free product with identification of units and $\psi:= \ast_{i \in I} \psi_i$ be the free product of linear functionals. $\varphi$ is defined by the following rule: if $X_k \in \mathcal{A}_{i_k}$ with $i_1 \neq \cdots \neq i_n$ and $\psi_{i_k}(X_k) = 0$ for all $1 \leq k \leq n$, then 
\begin{equation}\label{123}
\varphi (X_1 \cdots X_n) = \prod_{k = 1}^{n} \varphi_{i_k}(X_k).
\end{equation}
If $|I| =2$, we denote the c-free product as $(\varphi_1,\psi_1)\ast (\varphi_2,\psi_2) = (\varphi_1 {}_{\psi_1}\!\!\ast_{\psi_2}\!\varphi_2,\psi_1 \ast \psi_2)$, omitting the algebras for simplicity. \\ 
(2) Let $(\mathcal{A}, \varphi, \psi)$ be an algebraic probability space. Subalgebras $\{\mathcal{A}_i \}_{i \in I}$ are said to be c-free independent if they are free independent with respect to $\psi$ and satisfy the following property: 
\begin{equation}\label{1231}
\varphi (X_1 \cdots X_n) = \prod_{k = 1}^{n} \varphi(X_k)
\end{equation}
whenever $X_k \in \mathcal{A}_{i_k}$ with $i_1 \neq \cdots \neq i_n$ and $\psi(X_k) = 0$ for all $1 \leq k \leq n$. 
\end{defi}

Let $(\mathcal{A}, \varphi, \psi)$ be an algebraic probability space and $\comp[x]$ be the algebra generated by an indeterminate $x$. A distribution is a linear functional or sometimes a pair of linear functionals on the algebra $\comp[x]$. A distribution of $X \in \mathcal{A}$ is a single distribution $\mu_X$ or $\nu_X$, or sometimes a pair of distributions $(\mu_X, \nu_X)$ on $\comp[x]$ defined by $\mu_X(x^n) = \varphi(X^n)$ and $\nu_X(x^n) = \psi(X^n)$. Without mentioning explicitly, we hereafter assume that the symbols $\mu$ and $\nu$ denote distributions corresponding to $\varphi$ and $\psi$, respectively. If $X$ is unitary and selfadjoint in a $C^\ast$-algebra with a state, we can respectively identify $(\mu_X, \nu_X)$ with a pair of probability measures on $\tor$ and on $\real$. 

A multiplicative c-free convolution of probability measures has been introduced and investigated in \cite{P-W}. 
If $X$ and $Y$ are c-free independent, the distribution $(\mu_{XY},\nu_{XY})$ only depends on the distributions $(\mu_X, \nu_X)$ and $(\mu_Y, \nu_Y)$ on $\comp[x]$, not on a specific choice of an algebraic probability space or random variables. Therefore, 
we call $(\mu_{XY},\nu_{XY})$ a multiplicative c-free convolution of $(\mu_X,\nu_X)$ and $(\mu_Y,\nu_Y)$. We use the notation $(\mu_{XY},\nu_{XY}) = (\mu_X {}_{\nu_X}\!\!\boxtimes_{\nu_Y} \!\mu_Y, \nu_X \boxtimes \nu_Y) = (\mu_X, \nu_X) \boxtimes_c (\mu_Y, \nu_Y)$ for the multiplicative c-free convolution as a binary operation. Forgetting the random variables $X$ and $Y$, we can formulate the multiplicative c-free convolution $(\mu_1, \nu_1) \boxtimes_c (\mu_2, \nu_2)$ of two pairs of distributions $(\mu_i,\nu_i)$.  

Let $(\mu, \nu)$ be a pair of distributions on $\comp[x]$. 
We consider generating functions in the sense of formal power series. However, if distributions are bounded in such a way as $|\mu(x^n)| \leq A^n$ for a constant $A > 0$, then all the generating functions make sense as analytic functions. 
Let $G_\mu$ denote the Cauchy transform: $G_\mu(z)= \sum_{n = 0} ^\infty \frac{\mu(x^n)}{z^{n+1}}$.  
We also define $\eta_\mu(z):= 1 - \frac{z}{G_\mu (\frac{1}{z})}$ which plays important roles in descriptions of multiplicative convolutions in the free, Boolean, monotone cases. The $R$-transform and c-free $R$-transform are then defined from the functional relations 
\begin{align}
&\frac{1}{G _{\nu}(z)} = z - R_\nu (G_{\nu}(z)), \label{rel1}\\
&\frac{1}{G_{\mu}(z)} = z - R_{(\mu, \nu)}(G_{\nu}(z)). \label{rel2}  
\end{align}
If we introduce $\widetilde{R}_{(\mu, \nu)}(z):=z R_{(\mu, \nu)}(z)$ and $\widetilde{R}_{\mu}(z):=z R_{\mu}(z)$, the relations (\ref{rel1}) and (\ref{rel2})  are written as follows: 
\begin{align}
&\widetilde{R}_\nu \Big(\frac{z}{1 - \eta_\nu(z)}\Big) = \frac{\eta_\nu(z)}{1 - \eta_\nu(z)}, \label{eq01}\\
&\widetilde{R}_{(\mu, \nu)}\Big(\frac{z}{1 - \eta_\nu(z)}\Big) = \frac{\eta_\mu(z)}{1 - \eta_\nu(z)}. \label{eq02}   
\end{align}
These relations are more convenient than (\ref{rel1}) and (\ref{rel2}) in this paper.

 We define a c-free $T$-transform $T_{(\mu, \nu)}(z):= \frac{\widetilde{R}_{(\mu, \nu)}(\widetilde{R}^{-1} _\nu (z))}{\widetilde{R}^{-1} _\nu (z)}$ 
and a $T$-transform $T_\nu(z):= \frac{z}{\widetilde{R}^{-1}_\nu (z)}$.  
In the paper \cite{P-W}, Wang and Popa proved that the multiplicative c-free convolution is characterized by 
\begin{align}
&T_{(\mu_1, \nu_1)\boxtimes_c (\mu_2,\nu_2)}(z) = T_{(\mu_1, \nu_1)}(z) T_{(\mu_2, \nu_2)}(z), \\ 
&T_{\nu_1 \boxtimes \nu_2}(z) = T_{\nu_1}(z) T_{\nu_2}(z). 
\end{align}

\section{Observations on conditionally free independence and other notions of independence} \label{sec1}
We unify several multiplicative convolutions in the literature in terms of c-free convolutions. 
We denote the free, Boolean, monotone and orthogonal products by $\ast, \diamond$, $\rhd$ and $\angle$, respectively;  
for instance, the reader is referred to \cite{Len1,Mur5} for their definitions.  
We consider algebraic probability spaces $(\mathcal{A}_1, \varphi_1, \psi_1)$ and $(\mathcal{A}_2, \varphi_2, \psi_2)$. 
If $\mathcal{A}_i$ admits an algebra homomorphism $\delta^i: \mathcal{A}_i \to \comp$, 
then $\mathcal{A}_i$ has a decomposition 
\begin{equation}\label{assumption11}
\mathcal{A}_i = \comp 1 \oplus \mathcal{A}_i ^0
\end{equation} 
with $\mathcal{A}_i ^0:=\text{Ker}\,\delta^i$ $(i = 1, 2)$. In this case we have the following relations. 
\begin{gather}
(\varphi, \varphi) \ast (\psi, \psi) = (\varphi \ast \psi, \varphi \ast \psi) \text{~on~} \mathcal{A}_1  \ast \mathcal{A}_2, \label{free2} \\
(\varphi, \delta^1) \ast (\psi, \delta^2) = (\varphi \diamond \psi, \delta^1 \ast \delta^2) \text{~on~} \mathcal{A}_1 ^0  \ast \mathcal{A}_2 ^0,  \label{boolean2}\\
(\varphi, \delta^1) \ast (\psi, \psi) = (\varphi \rhd \psi, \psi) \text{~on~} \mathcal{A}_1 ^0 \ast \mathcal{A}_2, \label{monotone2} \\
(\varphi, \delta^1) \ast (\delta^2, \psi) = (\varphi \angle \psi, \psi) \text{~on~} \mathcal{A}_1 ^0 \ast \mathcal{A}_2. \label{ortho5} 
\end{gather}
The relations (\ref{free2}), (\ref{boolean2}) were found in \cite{BS, BLS}, the relation (\ref{monotone2}) in \cite{Fra} and the relation (\ref{ortho5}) in \cite{Has3}.

In the special case where $\mathcal{A}_1 = \comp[x_1]$ and $\mathcal{A}_2 = \comp[x_2]$, we can define a linear functional $\delta^j _c$ ($j = 1, 2,~c \in \comp$) by   
\[
\delta^j _c (x_j ^n):= c^n.  
\]
We note that $\delta^j = \delta^j _0$ holds. Then we obtain the following results. 
\begin{prop}\label{prop1} We denote by $(\varphi_1 {}_{\psi_1}\!\!\ast_{\psi_2}\!\varphi_2, \psi_1 \ast \psi_2)$ the c-free product $(\varphi_1, \psi_1) \ast (\varphi_2, \psi_2)$. 
Let $\mathcal{A}_1 = \comp[x_1]$ and $\mathcal{A}_2 = \comp[x_2]$. Let $(\varphi_j, \psi_j)$ be a pair of linear functionals on $\mathcal{A}_j$. 
\begin{itemize}
\item[(1)] $x_1 - c_1$ and $x_2 - c_2$ are Boolean independent in $(\mathcal{A}_1 \ast \mathcal{A}_2, \varphi_1 {}_{\delta^1_{c_1}}\!\!\ast_{\delta^2_{c_2}}\!\varphi_2)$. 

\item[(2)] $x_1 - c$ and $x_2$ are monotone independent in $(\mathcal{A}_1 \ast \mathcal{A}_2, \varphi_1 {}_{\delta^1_c}\!\!\ast_{\varphi_2}\!\varphi_2)$. 

\item[(3)] $x_1 - c_1$ and $x_2 - c_2$ are orthogonal independent in $(\mathcal{A}_1 \ast \mathcal{A}_2, \varphi_1 {}_{\delta^1_{c_1}}\!\!\ast_{\varphi_2}\!\delta^2_{c_2},~ \delta^1 _{c_1} \rhd \varphi_2)$. 

\item[(4)] $x_1 - c$ and $x_2$ are c-monotone independent in $(\mathcal{A}_1 \ast \mathcal{A}_2, \varphi_1 {}_{\delta^1 _c}\!\!\ast_{\psi_2}\!\varphi_2, \psi_1 \rhd \psi_2)$. 
\end{itemize}
\end{prop}
\begin{proof}
If we use the property $\delta^j_{c}((x_j -c)^n) = 0$ ($n \geq 1$), (1), (2) and (3) are equivalent to (\ref{boolean2}), (\ref{monotone2}) and (\ref{ortho5}), respectively. (4) follows from Theorem 3.6 of \cite{Has3}. 
 \end{proof}
 
In the commutative algebra $\comp[x]$, a linear functional is equivalent to a distribution of $x$. Then we can formulate multiplicative convolutions of distributions on $\comp[x]$, as explained in the c-free case. We re-write the equalities (\ref{boolean2})-(\ref{ortho5}) in terms of the multiplicative c-free convolution of distributions by using Proposition \ref{prop1}:  
\begin{gather}
(\mu, \delta_1) \boxtimes_c (\nu, \delta_1) = (\mu \utimes \nu, \delta_1), \label{boolean1}\\
(\mu, \delta_1) \boxtimes_c (\nu, \nu) = (\mu \dismrhd \nu, \nu), \label{monotone1} \\
(\mu, \delta_1) \boxtimes_c (\delta_1, \nu) = (\mu \angle \nu, \nu), \label{ortho1} 
\end{gather}
where $\autimes$, $\dismrhd$ and $\angle$ respectively denote the Boolean \cite{Fra2}, monotone \cite{Ber1} and orthogonal convolutions \cite{Len2} of distributions. 
We note that a symbol for a product of states sometimes differs from that for a convolution of distributions.

In papers \cite{Ber1, Ber2}, Bercovici has defined other multiplicative convolutions by supposing $X- \varphi(X)$ and $Y - \varphi(Y)$ are independent. We denote these convolutions of distributions by $\dismrhd_0$ and $\putimes_{\!0}$ in the monotone and Boolean cases, respectively. These are associative and characterized by the relations  
\begin{equation}
\eta_{\mu \submrhd_0 \nu}(z) = \eta_\mu \Big(\frac{1}{m_1(\mu)}\eta_\nu(m_1(\mu)z)\Big),~~~~~\frac{\eta_{\mu \hutimes_{\,0} \nu}(z)}{z}= \frac{\eta_{\mu}(m_1(\nu)z)}{m_1(\nu)z}\frac{\eta_{\nu}(m_1(\mu)z)}{m_1(\mu)z}, 
\end{equation}
where $m_n(\mu)$ is the $n$th moment $\mu(x^n)$. 
From Proposition \ref{prop1}, these convolutions are also written in the c-free setting: 
\begin{gather}
(\mu, \delta_{m_1(\mu)}) \boxtimes_c (\nu, \delta_{m_1(\nu)}) = (\mu \utimes_{\!0} \nu, \delta_{m_1(\mu)m_1(\nu)}), \label{boolean3}\\
(\mu, \delta_{m_1(\mu)}) \boxtimes_c (\nu, \nu) = (\mu \mrhd_0 \nu, T_{m_1(\mu)}\nu), \label{monotone3} 
\end{gather}
where $T_c \mu$ is characterized by $\eta_{T_c \mu}(z) = c\eta_{\mu}(z)$. 
In addition, if we use a multiplicative c-monotone convolution which will be introduced in Definition \ref{defi222},  
\begin{equation}
(\mu, \delta_{m_1(\mu)}) \mrhd_c (\nu, \delta_{m_1(\nu)}) = ((S_{m_1(\nu)}\mu) \utimes \nu, \delta_{m_1(\mu)m_1(\nu)}), \label{monotone31}\\
\end{equation}
where $S_c \mu$ is defined by $\eta_{S_c \mu}(z) = \frac{\eta_\mu(cz)}{c}$. Therefore, the associative law of $\putimes_{\!0}$ is naturally understood in terms of $\boxtimes_c$; we can say that $\dismrhd_0$ is the multiplicative version of the Fermi convolution \cite{O1}. By contrast, the associative law of $\dismrhd_0$ cannot be understood in terms of $\boxtimes_c$ or $\dismrhd_c$. We do not treat this problem in this paper.

\section{Multiplicative conditionally monotone convolutions} \label{sec2}
A multiplicative convolution for c-monotone independence is defined as follows. 
\begin{defi}\label{defi222} Consider an algebraic probability space $(\mathcal{A}, \varphi, \psi)$. 
Let $X, Y$ be elements of $\mathcal{A}$ such that $X - 1$ and $Y$ are c-monotone independent (or equivalently, $X-1$ and $Y-1$ are c-monotone independent). Then a multiplicative c-monotone (or $\dismrhd_c$- for short) convolution is defined by the distribution of $XY$. 
\end{defi}
The reason why we consider not $X$ but $X-1$ can be partially understood from Proposition \ref{prop1} and the relations (\ref{boolean1})-(\ref{ortho1}). However, it is expected that the reason is more clarified in future researches.  

By definition, the left distribution $\mu_{XY}$ only depends on $\mu_X, \mu_Y, \nu_Y$. The right distribution $\nu_{XY}$ is the  multiplicative monotone convolution. Therefore, we denote them as $(\mu_{XY},\nu_{XY}) = (\mu_X \dismrhd_{\nu_Y}\mu_Y, \nu_X \dismrhd \nu_Y) = (\mu_X,\nu_X) \dismrhd_c (\mu_Y,\nu_Y)$. As is the case for other convolutions,  we can only consider distributions on $\comp[x]$, forgetting the random variables $X$ and $Y$.

An immediate consequence of Proposition \ref{prop1} is a connection to the c-free convolution.  
\begin{prop}\label{cor312} 
For distributions $\mu_1, \mu_2, \nu_2$, we have $\mu_1 \dismrhd_{\nu_2}\mu_2 = \mu_1 {}_{\delta_1}\!\!\boxtimes_{\nu_2} \mu_2$.  
\end{prop}

Proposition \ref{cor312} enables us to characterize the c-monotone convolution in terms of analytic functions used in the c-free case. The result, however, is not trivial.  
\begin{thm} \label{thm01} 
For distributions $\mu_1,\nu_2,\mu_2,\nu_2$ on $\comp[x]$, we have  
\begin{align}
\eta_{\mu_1 \submrhd_{\nu_2} \mu_2 }(z) = \frac{\eta_{\mu_2}(z)}{\eta_{\nu_2}(z)}\eta_{\mu_1}(\eta_{\nu_2}(z)), \label{eq422}\\
\eta_{\nu_1 \submrhd \nu_2}(z) = \eta_{\nu_1}(\eta_{\nu_2}(z)). \label{eq42}
\end{align}
If $\eta_{\nu_2}=0$, for instance if $\nu_2$ corresponds to the normalized Haar measure $\omega$ on $\tor$, (\ref{eq422}) is understood to be $\eta_{\mu_1 \submrhd_{\nu_2} \mu_2}(z) = \frac{d}{dz}\eta_{\mu_1}(0)\eta_{\mu_2}(z) =\mu_1(x)\eta_{\mu_2}(z) $.  

If distributions are arising from probability measures on the unit circle, the above relations are valid as analytic maps for $|z| < 1$.   
\end{thm}
\begin{proof}
First we assume that the mean $\nu_2(x)$ is non-zero since we use the inverse function of $\widetilde{R}_{\nu_2}$ whose coefficient of $z$ is equal to $\nu_2(x)$.   
 (\ref{eq42}) was proved in \cite{Ber1}. 
We note that $T_{(\mu, \delta_1)} (z) = \frac{1 +z}{z}\eta_\mu \big( \frac{z}{1+z} \big)$. Then 
we have 
\begin{equation}\label{eq43}
\begin{split}
T_{(\mu_1 \submrhd_{\nu_2}\mu_2, \nu_2)}(z) 
           &= T_{(\mu_1, \delta_1)\boxtimes_c(\mu_2, \nu_2)}(z) \\
           &=  T_{(\mu_1, \delta_1)}(z)T_{(\mu_2, \nu_2)}(z) \\
           &= \frac{1 +z}{z}\eta_{\mu_1} \Big( \frac{z}{1+z} \Big) \frac{\widetilde{R}_{(\mu_2, \nu_2)}(\widetilde{R}^{-1} _{\nu_2} (z))}{\widetilde{R}^{-1} _{\nu_2} (z)}.   
\end{split}
\end{equation}
Therefore, we have 
\begin{equation}\label{eq431}
\widetilde{R}_{(\mu_1 \submrhd_{\nu_2}\mu_2, \nu_2)}(\widetilde{R}^{-1} _{\nu_2} (z)) = \frac{1 +z}{z}\eta_{\mu_1} \Big( \frac{z}{1+z} \Big) \widetilde{R}_{(\mu_2, \nu_2)}(\widetilde{R}^{-1} _{\nu_2} (z)).
\end{equation}
We define $w$ by the relation $\widetilde{R}^{-1} _{\nu_2} (z) = \frac{w}{1 - \eta_{\nu_2}(w)}$. This is equivalent to 
$z = \widetilde{R}_{\nu_2}\big( \frac{w}{1 - \eta_{\nu_2}(w)} \big) = \frac{\eta_{\nu_2}(w)}{1 - \eta_{\nu_2}(w)}$. 
Then we have 
\begin{equation}\label{eq44}
\frac{z}{1 + z} = \eta_{\nu_2}(w). 
\end{equation}
Combining the equalities (\ref{eq02}), (\ref{eq431}), (\ref{eq44}), we obtain the conclusion.  

Second, we consider the case $\nu_2(x)=0$. We note that the moments of $\mu_1 \dismrhd_{\nu_2} \mu_2$ can be expressed in terms of  sums and products of moments of $\mu_1, \mu_2, \nu_2$. Therefore we can approximate $\nu_2$ by a sequence $\nu_2^{(n)}$, $n=1,2,3,\cdots$, so that $\nu_2^{(n)}(x) \neq 0$ and $\nu_2^{(n)}$ converges to $\nu_2$ in the sense of moments. 
We let $n$ tend to infinity and then the conclusion follows. This proof includes the case $\eta_{\nu_2} =0$. 
\end{proof}

This characterization includes $\dismrhd$, $\putimes$ and the multiplicative orthogonal convolution $\angle$: these convolutions have been characterized in \cite{Ber1, Fra, Len2} as 
\begin{align}
\eta_{\mu \hutimes \nu}(z) = \frac{\eta_{\mu}(z)\eta_{\nu}(z)}{z}, \\
\eta_{\mu \angle \nu}(z) = \frac{z\eta_{\mu}(\eta_{\nu}(z))}{\eta_{\nu}(z)}, \\ 
\eta_{\mu \submrhd \nu} = \eta_{\mu} \circ \eta_{\nu}. 
\end{align} 

 In terms of c-monotone convolutions, the monotone, Boolean and orthogonal convolutions appear as follows. 
\begin{gather}
(\mu, \nu) \mrhd_c(\mu, \nu) = (\mu \mrhd \nu, \mu \mrhd \nu), \label{monotone6} \\
(\mu, \delta_1) \mrhd_c (\nu, \delta_1) = (\mu \putimes \nu, \delta_1), \label{boolean6} \\ 
(\mu, \lambda) \mrhd_c (\delta_1, \nu) = (\mu \angle \nu, \lambda \mrhd \nu). 
\end{gather}

\section{Convolution semigroups for multiplicative convolutions}\label{sec3}
From this section, we consider distributions coming from probability measures.  
We respectively denote by $\tor$ and $\disc$ the unit circle $\{z \in \comp:|z|=1 \}$ and the unit disc $\{z \in \comp: |z| < 1 \}$. 
Moreover, let $\pa(\tor)$ and $\pa(\real_+)$ be the sets of probability measures on $\tor$ and $\real_+$, respectively.  
It is known that the multiplicative monotone convolution and orthogonal convolution preserve $\pa(\real_+)$ \cite{Ber1, Len2}.  
The multiplicative Boolean convolution, however, does not preserve $\pa(\real_+)$, and hence, the multiplicative c-monotone convolution does not, either (see \cite{Ber2} and also \cite{Fra3}). We do not investigate into this problem in this paper and we focus on probability measures on $\pa(\tor)$ from now on. 

The following characterization is known (see \cite{Ber-Bel}). 
\begin{prop}\label{prop3}
Let $\eta: \disc \to \comp$ an analytic function. The following conditions are equivalent. \\
(1) There exists a probability measure $\mu \in \pa(\tor)$ such that $\eta = \eta_\mu$. \\
(2) $\eta(0) = 0$ and $|\eta(z)| < 1$ for all $z \in \disc$.  \\
(3) $|\eta(z)| \leq |z|$ for all $z \in \disc$. 
\end{prop}
We know that $(\mu_1, \nu_1) \dismrhd_c (\mu_2, \nu_2) \in \pa(\tor) \times \pa(\tor)$ if $(\mu_1, \nu_1), (\mu_2, \nu_2) \in \pa(\tor) \times \pa(\tor)$, since the product of unitary elements is again a unitary. We can also prove this property from Proposition \ref{prop3} easily. 

The following points are useful to understand this paper. The convolution for right components is just a monotone convolution whose properties have been studied in details in the literature (see \cite{Ber1,Fra} and also \cite{Be-Po}). 
We often use such results to prove properties of left components. However, sometimes a proof for left components 
essentially includes a proof for right components 
if we set the probability measures of the left and right components equal. We have met such a situation in Theorem \ref{thm01}: the proof of Eq. (\ref{eq422}) actually generalizes Eq. (\ref{eq42}).

We prove a correspondence between a $\dismrhd_c$-convolution semigroup and a pair of vector fields. We consider a $\dismrhd_c$-convolution semigroup $\{(\mu_t, \nu_t) \}_{ t \geq 0}$ with $(\mu_0, \nu_0) = (\delta_1, \delta_1)$. If we define $F_t(z):= \log \eta_{\mu_t}(e^z)$ and $H_t(z):= \log \eta_{\nu_t}(e^z)$ in a suitable domain, we have the relations $F_{t + s}(z) = F_t (F_{s}(z))$ and $H_{t+s}(z) = H_t(F_s(z)) - F_s(z) + H_s(z)$. 
Interestingly these relations coincide with the additive c-monotone convolution case. However, we need to restrict the domain to define the logarithm and Muraki's method in \cite{Mur3} cannot be applied. 
We give a proof based on Berkson and Porta's result on composition semigroups \cite{Be-Po}.   

\begin{thm} \label{thm3} Let $U$ be an open set in $\comp$. Let $\{F_t(z) \}_{t \geq 0}$ and $\{H_t(z) \}_{t \geq 0}$ be families of analytic maps $F_t: U \to U$, $H_t: U \to U$ satisfying $F_0(z) = z$, $H_0(z) = z$,  $F_{t + s}(z) = F_t (F_{s}(z))$ and $H_{t+s}(z) = H_t(F_s(z)) - F_s(z) + H_s(z)$.  
We assume that $(t, z) \mapsto F_t(z)$ and $(t, z) \mapsto H_t(z)$ are both continuous on $[0,\infty) \times U$. 
Then there exist analytic vector fields $A_1$ and $A_2$ in $U$ such that 
\begin{align}
\frac{d}{dt} H_t(z) = A_1(F_t(z)), \label{diff1} \\
\frac{d}{dt} F_t(z) = A_2(F_t(z)),  \label{diff2}
\end{align}
for $z \in U$ and $t \in [0, \infty)$. 
\end{thm}
\begin{proof} The fact for $F_t$ and $A_2$ is known in \cite{Be-Po}. We prove the claim for $H_t$ and $A_1$.  
For any compact convex set $K \subset U$, there exists $\alpha > 0$ such that the convex hull of $\cup \{ F_t(K): t \in [0, \alpha] \}$ is a compact set in $U$. 
We denote the convex hull by $\widetilde{K}$. In this proof, we always use $C$ to mean that there exists some constant (dependent on $K$). 
Berkson and Porta have proved that 
\begin{equation}
|F_t(z) - z| \leq Ct^{\frac{2}{3}},~~z \in K,~~ t \in [0, 1].
\end{equation}
We apply the same method to $H_t$. A key equality is the following: 
\begin{equation}\label{eq32}
\begin{split}
H_{2t}(z) - 2H_t(z) + z &= \int_z ^{F_t(z)} \frac{d}{dw}(H_t(w) - w)dw \\ 
                        &=  \int_z ^{F_t(z)} dw \frac{1}{2\pi i}\int_\Gamma \frac{H_t(\zeta) - \zeta}{(\zeta - w)^2}d\zeta, 
\end{split}
\end{equation}
where $\Gamma$ is a closed curve around $\widetilde{K}$. The path for the integration with respect to $w$ is the line segment from $z$ to $F_t(z)$. 
By simple estimation we obtain 
\begin{equation}
|H_{2t}(z) - 2H_t(z) + z| \leq C|F_t(z) - z| \leq Ct^{\frac{2}{3}}
\end{equation}
for $t \in [0, 1]$ and $z \in K$. 
Then we have 
\begin{equation}
|H_{t}(z) -  z| \leq \frac{1}{2}|H_{2t}(z) - z| + Ct^{\frac{2}{3}} ~~\text{for}~~ t \in [0, 1],~~z \in K.
\end{equation}
By iteration, we have 
\begin{equation}\label{eq310}
\begin{split}
|H_{2^{-n}t}(z) - z| &\leq 2^{-n}|H_t(z) - z| + C t^{\frac{2}{3}}2^{-\frac{2}{3}n}\sum_{k = 1} ^n 2^{-\frac{k}{3}} \\ 
                     &\leq  2^{-n}|H_t(z) - z| + C t^{\frac{2}{3}}2^{-\frac{2}{3}n}
\end{split}
\end{equation}
for $t \in [0, 1]$ and $z \in K$, $n \geq 0$. Since $\{2^{-n}t: n \geq 0,~~t \in [\frac{1}{2}, 1] \} = [0,1]$, the estimate 
\begin{equation}\label{eq31}
|H_t(z) - z| \leq Ct^{\frac{2}{3}},~~z \in K,~~ t \in [0, 1]
\end{equation} 
follows. More precisely, let $s= 2^{-n}t$ for $n \geq 0$ and $\frac{1}{2} \leq t \leq 1$. 
Then (\ref{eq310}) can be written as $|H_s(z)-z| \leq \frac{s}{t}|H_t(z)-z| + Cs^{\frac{2}{3}}$. Since $|\frac{H_t(z)-z}{t}|$ is bounded for 
$t \in [\frac{1}{2},1]$ and $z \in K$, we have $|H_s(z) -z| \leq Cs + Cs^{\frac{2}{3}}$. This inequality holds for all $s \in [0, 1]$. Finally, we 
note that $s \leq s^{\frac{2}{3}}$ for $s \in [0,1]$.

The remaining discussion is the same as the original paper. We do not repeat the argument. 
\end{proof} 
%
%

To prove the main theorem in this section, we need the following fact.  
We note that the proof does not depend on the semigroup property. 
\begin{prop}\label{prop10}
Let $\{\phi_t \}_{t \in I}$ be a family of analytic maps on $\disc$ parametrized by $t \in I$, where $I$ is an interval. We assume that the map $t \mapsto \phi_t(z)$ 
is continuous for each $z \in \disc$. Then the map $\phi: [0, \infty) \times \disc \to \disc$ defined by $\phi(t,z) = \phi_t(z)$ 
is continuous.  
\end{prop}
\begin{proof}
Let $B_R:=\{z \in \disc: |z| < R \}$ for $R < 1$. By Cauchy's integral formula, we have 
\[
\phi_t(z) = \frac{1}{2\pi i}\int_{\partial B_R} \frac{\phi_t(w)}{w -z}dw
\]
for $z \in B_R$ and $t \in I$.  Let $(t_n, z_n)$ be a sequence converging to $(t, z)$. By Lebesgue's bounded convergence theorem, 
$\phi_{t_n}(z_n) \to \phi_t(z)$ as $n \to \infty$ since $|\phi_{t}(w)| < 1$.  
\end{proof}

Now consider a weakly continuous $\dismrhd_c$-convolution semigroup $\{(\mu_t, \nu_t) \}_{ t \geq 0}$ with $(\mu_0, \nu_0) = (\delta_1, \delta_1)$. From the weak continuity, $\eta_{\mu_t}$ and $\eta_{\nu_t}$ are continuous as functions of $t$ for each $z$. Moreover, $\eta_{\mu_t}(z)$ and $\eta_{\nu_t}(z)$ are continuous in $C([0,\infty) \times \disc)$ from Proposition \ref{prop10}. 
We take a compact disc $D \subset \disc$, $0 \notin D$. Without loss of generality we assume that $D \subset \com+$. There exists $\alpha$ such that $\cup \{\eta_{\mu_t}(D), \eta_{\nu_t}(D); t \in [0, \alpha] \} \subset \com+$ and then we may define $\log \eta_{\mu_t}(z)$ and $\log \eta_{\nu_t}(z)$ for $z \in D,~~t \in [0, \alpha]$. We define $F_t(z):= \log \eta_{\nu_t}(e^z)$ and $H_t(z):= \log \eta_{\mu_t}(e^z)$ and their domain $E:= \log(D)$. The images $H_t(E)$ and $F_t(E)$ may not be included in $E$, but we can use the technique of Theorem \ref{thm3} for small $t$ and obtain the differentiability of $H_t$. 

We quote the following theorem; the reader is referred to \cite{Akh}. 
\begin{thm}\label{thm6}
Let $f: \disc \to \{z \in \comp: \re z \geq 0 \}$ be an analytic function. Then $f$ can be represented as 
\[
f(z) = ib + \int_{-\pi} ^\pi \frac{e^{i\theta} + z}{ e^{i\theta} - z}\rho(d\theta), 
\]
where $b \in \real$ and $\rho$ is a positive finite Borel measure. Then $b = \im f(0)$ and 
$\rho([\alpha, \beta]) = \lim_{r \nearrow 1}\frac{1}{2\pi}\int_\alpha ^\beta \re f(re^{i\theta})d\theta$ 
for all continuity points $\alpha, \beta$ of $\rho$. 
\end{thm}
\begin{thm}\label{thm4}
Let $\{ (\mu_t, \nu_t) \}_{t \geq 0}$ be a weakly continuous $\dismrhd_c$-convolution semigroup with $(\mu_0, \nu_0) = (\delta_1, \delta_1)$. Then there exist analytic functions $B_1$, $B_2: \disc \to \comp$ satisfying  $\re B_1, \re B_2 \leq 0$ such that 
\begin{align}
\frac{d}{dt} \eta_{\mu_t}(z) = \eta_{\mu_t}(z) B_1(\eta_{\nu_t}(z)), \label{eq51} \\ 
\frac{d}{dt} \eta_{\nu_t}(z) = \eta_{\nu_t}(z) B_2(\eta_{\nu_t}(z)). \label{eq52}
\end{align}  

Conversely, if two analytic functions $B_1$, $B_2: \disc \to \comp$ are given satisfying $\re B_1, \re B_2 \leq 0$, there corresponds a weakly 
continuous $\dismrhd_c$-convolution semigroup $\{ (\mu_t, \nu_t) \}_{t \geq 0}$  with $(\mu_0, \nu_0) = (\delta_1, \delta_1)$ defined by (\ref{eq51}) and (\ref{eq52}). 

The vector fields are written in the Herglotz-Riesz formula 
\begin{equation}\label{levy}
B_j(z) = i\gamma_j + \int_{\tor}\frac{z + \zeta}{z - \zeta}\tau_j(d\zeta),~~ j = 1,2,
\end{equation}  
where $\gamma_j$ is a real number and $\tau_j$ is a positive finite measure. This formula is the analogue of the L\'{e}vy-Khintchine formula in probability theory.  
\end{thm}
\begin{proof}
For the second component, the claim is identical to the monotone case in \cite{Ber1} and we only explain the first component. The first component is similar. The existence of the vector fields is the consequence of Theorem \ref{thm3}, as explained in the above.  
By Proposition \ref{prop3}, $\eta_{\mu_t}$ satisfies
\begin{align}
|\eta_{\mu_{t + s}}(z)| &= \Big|\frac{\eta_{\mu_t}(\eta_{\nu_s}(z))}{\eta_{\nu_s}(z)}\Big||\eta_{\mu_s}(z)| \\ \notag
                  & \leq |\eta_{\mu_s}(z)|, \notag
\end{align}
which implies that $|\eta_{\mu_t}(z)|$ is a non-increasing function of $t$. By the way, 
(\ref{eq51}) implies that  
\begin{align}
\frac{d}{dt}|\eta_{\mu_t}(z)|^2 = 2|\eta_{\mu_t}(z)|^2\re B_1(\eta_{\nu_t}(z)).
\end{align}
Indeed, $\frac{d}{dt}|\eta_{\mu_t}(z)|^2 = (\frac{d}{dt}\eta_{\mu_t}(z)) \overline{\eta_{\mu_t}(z)} +\eta_{\mu_t}(z) \overline{ \frac{d}{dt}\eta_{\mu_t}(z)} = |\eta_{\mu_t}(z)|^2B_1(\eta_{\nu_t}(z)) + |\eta_{\mu_t}(z)|^2 \overline{B_1(\eta_{\nu_t}(z))}$. 
Therefore, $B_1$ needs to satisfy $\re B_1 \leq 0$. 

In the converse statement, the existence of $\eta_{\nu_t}$ is a consequence of the result in \cite{Ber1}. Therefore, we only needs to prove the existence of $\eta_{\mu_t}$. If a given vector field $B_1$ satisfies $\re B_1 \leq 0$, we can define $\kappa_{t}$ by 
\begin{equation}
\kappa_{t}(z) = z\exp\Big(\int_0 ^t B_1(\eta_{\nu_s}(z))ds\Big). 
\end{equation}
We can prove the functional equation $\kappa_{t + s}(z) = \frac{\kappa_{t}(\eta_{\nu_s}(z))}{\eta_{\nu_s}(z)}\kappa_{s}(z)$ as follows.  
Let $f_t(z):=\kappa_{t+s}(z)$ and $g_t(z):=\frac{\kappa_{t}(\eta_{\nu_s}(z))}{\eta_{\nu_s}(z)}\kappa_{s}(z)$ for a fixed $s \geq 0$. 
The differential equations for $f_t(z)$ and $g_t(z)$ are 
\[
\frac{d}{dt}f_t(z) = f_t(z)B_1(\eta_{\nu_{t+s}}(z))
\]
 and 
\[
\begin{split}
\frac{d}{dt}g_t(z) &= \kappa_{t}(\eta_{\nu_s}(z))B_1(\eta_{\nu_t}(\eta_{\nu_s}(z)))\frac{\kappa_{s}(z)}{\eta_{\nu_s}(z)} \\
             &= g_t(z)B_1(\eta_{\nu_{t+s}}(z)). 
\end{split}             
\] 
These two equations imply that $\frac{1}{f_t(z)}\frac{d}{dt}f_t(z) = \frac{1}{g_t(z)}\frac{d}{dt}g_t(z)$, and therefore $f_t(z) =g_t(z)$.   

$\re B_1 \leq 0$ implies that $|\kappa_t(z)|$ is non-increasing, and hence, $|\kappa_t(z)| \leq |\kappa_0(z)| \leq |z|$. By Proposition \ref{prop3}, there exists $\mu_t \in \pa(\tor)$ such that $\kappa_t = \eta_{\mu_t}$.  In conclusion, $(\mu_t, \nu_t)$ forms a c-monotone convolution semigroup. 
\end{proof}

\section{Infinitely divisible distributions}\label{sec4} 
\subsection{Embedding of  an infinitely divisible distribution to a convolution semigroup}
Infinitely divisible distributions form an important class of probability measures in probability theory. It is well known that 
an infinitely divisible distribution can be embedded into a continuous convolution semigroup. 
We establish the analogy for the multiplicative $c$-monotone convolution.  
We start from the definition of infinite divisibility. 
\begin{defi}
$(\mu, \nu) \in \pa(\tor) \times \pa(\tor)$ is said to be $\dismrhd_c$-infinitely divisible if and only if for any natural number $n \geq 2$, there exists $(\mu_n, \nu_n) \in \pa(\tor) \times \pa(\tor)$ such that $(\mu, \nu) = (\mu_n, \nu_n)^{\submrhd_c n}$. 
\end{defi}

From now on $\omega$ denotes the normalized Haar measure on $\tor$.
\begin{lem}\label{lem2}
Let $(\mu, \nu)$ be $\dismrhd_c$-infinitely divisible. \\
(1) If $\int_{\tor}\zeta d\mu(\zeta) = 0$, then $\mu = \omega$. \\
(2) If $\int_{\tor}\zeta d\nu(\zeta) = 0$, then $\nu = \omega$. 
\end{lem}
\begin{proof}
The fact (2) is known in \cite{Ber1}. First we prove the following fact. \\

$(\ast)$ Let $\lambda, \rho \in \pa(\tor)$. 
We define $a_k(\lambda)$ by $\eta_\lambda(z) = \sum_{k = 1}^\infty a_k(\lambda)z^k$ for $\lambda \in \pa(\tor)$. 
Let $(\lambda^n, \rho^n):= (\lambda, \rho)^{\submrhd_c n}$. If $\int_{\tor}\zeta d\lambda(\zeta) = 0$, we have  $a_1(\lambda^n) = \cdots = a_n(\lambda^n) = 0$. \\

We prove this by induction. If $n=1$, the statement is trivial. We assume that this property holds for $n=p$. Then 
\[
\begin{split}
\eta_{\lambda^{p+1}}(z)&=\eta_{\lambda \submrhd_{\rho^{p}}\lambda^p}(z) \\
                       &=\frac{\eta_{\lambda}(\eta_{\rho^p}(z))}{\eta_{\rho^p}(z)} \eta_{\lambda^p}(z) \\
 &= \Big( \sum_{k\geq 2}a_k(\lambda)\eta_{\rho^p}(z)^{k-1}\Big)\sum_{k\geq 1}a_k(\lambda^p)z^k.   
\end{split}
\]
Since $a_1(\lambda^p)= \cdots = a_{p}(\lambda^p)=0$ by assumption, the power of $z$ in $\eta_{\lambda^{p+1}}(z)$ starts from $p+2$. This implies that 
$a_1(\lambda^{p+1})=\cdots =a_{p+1}(\lambda^{p+1})=0$.

Let $\mu_n, \nu_n \in \pa(\tor)$ ($n \geq 2$) be probability measures such that $(\mu, \nu):= (\mu_n, \nu_n)^{\submrhd_c n}$. We observe first that $ 0 = a_1(\mu) = a_1(\mu_n)^n$, which implies $a_1(\mu_n) = 0$ for any $n \geq 2$. 
Then we can apply the above fact to conclude (1).  
\end{proof}
\begin{thm}\label{thm2121}
Let $\mu, \nu \in \pa(\tor)$ be probability measures. The following statements are equivalent. 
\begin{itemize}
\item[(1)] $(\mu, \nu)$ is $\dismrhd_c$-infinitely divisible with $\mu  \neq \omega, \nu \neq \omega$. 
\item[(2)] There exists a weakly continuous $\dismrhd_c$-convolution semigroup $\{(\mu_t, \nu_t) \}_{t \geq 0}$ with $(\mu_0, \nu_0) = (\delta_1, \delta_1)$ and $(\mu_1, \nu_1) = (\mu, \nu)$. 
\end{itemize}
\end{thm}

\begin{proof}
The proof is the same as the monotone case (see Theorem 4.4 in \cite{Ber1}) if we use Lemma \ref{lem2}, and we omit the proof. 
\end{proof}
\begin{rem}\label{rem3}
The convolution semigroup $\{(\mu_t,\nu_t)\}$ in the statement (2) is not unique as we will show in Subsection \ref{subsec5}. 
\end{rem}

The following properties are useful to understand the $\dismrhd$- and 
$\dismrhd_c$-convolutions. 

\begin{prop}\label{prop21}
Let $\nu$ be a delta measure at a point in $\tor$ and $\{\nu_t \}_{t \geq 0}$ be a $\dismrhd$-convolution semigroup as in Theorem \ref{thm2121} (2). 
Then, the associated function $B_2(z)$ satisfying  (\ref{eq52}) is a constant with value in $i\real$.  
\end{prop}
\begin{proof}
Let $\nu$ be the delta measure $\delta_{e^{i\alpha}}$ for an $\alpha \in \real$. By expanding $\eta_{\nu_t}(z) = \sum_{n=1}^\infty b_n(t)z^n$ and $B_2(z)=\sum_{n=1}^\infty s_nz^{n-1}$ in (\ref{eq52}), 
we have $\frac{d}{dt}b_1(t) = s_1b_1(t)$. The initial condition is $b_1(0) =1$, so that $b_1(t) = e^{s_1 t}$. Since $b_1(1) = e^{i\alpha}$, $s_1 = i\alpha + 2\pi i n$ for an integer $n$. 
In particular, $s_1 \in i\real$.  By the way, the L\'{e}vy-Khintchine formula (\ref{levy}) for $B_2$ implies that 
$\re s_1 = -\tau_2(\tor)$. Therefore, $\tau_2 = 0$, and the function $B_2$ is a constant in $i\real$. 
 \end{proof}
\begin{rem}In connection to Remark \ref{rem3}, $B_2(z) = i\alpha + 2 \pi i n$ ($n \in \mathbb{Z}$) generates the same probability measure $\delta_{e^{i\alpha}}$ at time $t=1$.  
The translation by  $2\pi i n$ however  does not preserve the probability measure at time 1 in generic cases. In the next subsection we will investigate this problem more. 
\end{rem}
\begin{prop}\label{prop31}
Let $\{(\mu_t, \nu_t) \}_{t \geq 0}$ be a weakly continuous $\dismrhd_c$-convolution semigroup with $(\mu_0, \nu_0) = (\delta_1, \delta_1)$. 

(1) If  $B_2 \equiv 2\pi i \frac{p}{q}$ for integers $p \neq 0$ and $q >0$ which cannot be divided by a common prime number, then 
$\eta_{\mu_1}(z) = z\exp(r_1+\sum_{n \in \nat \cap(q\nat)^c} r_{n + 1} \frac{e^{2\pi i np/q}-1}{2 \pi i np /q}z^{n})$ where $\nat =\{1,2,3,\cdots\}$. 
In particular, if $B_2 \equiv 2 \pi i k$ for a non-zero integer $k$, $\eta_{\mu_1}$ becomes $e^{r_1}z$, and the explicit density function of $\mu_1$ is the Poisson kernel shown in Example \ref{exa3}. 

(2) If $B_2 \equiv 0$, then $\eta_{\mu_1}(z) = ze^{B_1(z)}$. 
\end{prop}
\begin{proof}
If $B_2 \equiv 2 \pi i \frac{p}{q}$, $\eta_{\mu_1}$ is expressed as 
\begin{equation*}
\eta_{\mu_1}(z) = z\exp\Big(\int_0 ^1 B_1(\eta_{\nu_s}(z))ds\Big) = z \exp\Big(\int_0 ^1 B_1(e^{2 \pi i  \frac{sp}{q}}z)ds\Big)   
\end{equation*} 
by (\ref{eq51}). 
We expand the function $B_1$ as $B_1(z) = \sum_{n=1}^\infty r_n z^{n-1}$ and the integral  becomes 
\[
\begin{split}
\int_0 ^1 B_1(e^{2 \pi i \frac{sp}{q}}z)ds &= \sum_{n=0}^\infty z^{n} r_{n+1}\int_0 ^1 e^{2\pi i \frac{spn}{q}}ds.   
\end{split}
\]
The integral in the RHS vanishes if and only if $n$ can be divided by $q$. If $n$ cannot be divided by $q$, the integral $\int_0 ^1 e^{2\pi i \frac{spn}{q}}ds$ becomes $\frac{e^{2\pi i np/q}-1}{2\pi i np/q}$. 

If $B_2 \equiv 0$, then we can prove the claim by the same method. 
\end{proof}


\subsection{On convolution semigroups which have the same distribution at time one} \label{subsec5}
We denote by $ID(\putimes, \tor)$ the set of all $\putimes$-infinitely divisible distributions and define 
$ID(\putimes, \tor)_0 := ID(\putimes, \tor) \backslash \{\omega \}$. 
Franz proved in \cite{Fra2} that a probability measure $\mu \in \pa(\tor)$ belongs to $ID(\putimes, \tor)_0$ 
if and only if $\frac{\eta_\mu(z)}{z}$ (defined by $\eta_\mu '(0)$ at the origin) 
does not have a zero point in $\disc$. This condition is equivalent to 
the condition that there exists an analytic map 
$B: \disc \to \{z \in \comp: \re z \leq 0 \}$ such that $\eta_\mu(z) = z e^{B(z)}$. 
The above two conditions are also equivalent to the condition that $\mu$ can be embedded into a convolution semigroup $\{\mu_t \}_{t \geq 0}$. 
The relation between $B(z)$ and $\mu_t$ is $\eta_{\mu_t}(z)= ze^{tB(z)}$. We can understand this relation as a special case of (\ref{eq52}) where $\nu_t$ are all equal to $\delta_1$. 
Now there is a problem which does not arises in the additive convolution: the function $B$ is not unique. 
The non-uniqueness is however only due to the transformations
\begin{equation}\label{eq611}
B \mapsto B + 2 \pi i n \text{~for~} n \in \mathbb{Z}. 
\end{equation}

We consider the problem of uniqueness in the monotone and c-monotone cases. 
We follow the notation in Theorem \ref{thm4}. 
Let $\{(\mu_t,\nu_t) \}_{t \geq 0}$ and $\{(\widetilde{\mu}_t,\widetilde{\nu}_t) \}_{t \geq 0}$ be weakly continuous $\dismrhd_c$-convolution 
semigroups satisfying $(\mu_0,\nu_0) = (\widetilde{\mu}_0, \widetilde{\nu}_0) = (\delta_1,\delta_1)$ and 
$(\mu_1,\nu_1) = (\widetilde{\mu}_1,\widetilde{\nu}_1)$. 
The vector fields for $\{(\widetilde{\mu}_t,\widetilde{\nu}_t) \}_{t \geq 0}$ is denoted by $(\widetilde{B}_1,\widetilde{B_2})$. 
We assume that all $\mu_1,\nu_1,\widetilde{\mu}_1,\widetilde{\nu}_1$ are different from the normalized Haar measure.  
In addition we expand the four vector fields as  
\begin{equation}
\begin{split}
&B_1(z) = \sum_{n=1}^\infty r_n z^{n-1},~~~~~~~B_2(z)= \sum_{n=1}^\infty s_n z^{n-1}, \\
&\widetilde{B}_1(z) = \sum_{n=1}^\infty \widetilde{r}_n z^{n-1},~~~~~~~\widetilde{B}_2(z)= \sum_{n=1}^\infty \widetilde{s}_n z^{n-1}. 
\end{split}
\end{equation}
Also we expand $\eta_{\mu_t}$, $\eta_{\nu_t}$,  $\eta_{\widetilde{\mu}_t}$ and $\eta_{\widetilde{\nu}_t}$ as
\begin{equation}
\begin{split}
&\eta_{\mu_t}(z) = \sum_{n=1}^\infty a_n(t)z^n,~~~~~~~\eta_{\nu_t}(z) = \sum_{n=1}^\infty b_n(t)z^n, \\
&\eta_{\widetilde{\mu}_t}(z) = \sum_{n=1}^\infty \widetilde{a}_n(t)z^n,~~~~~~~\eta_{\widetilde{\nu}_t}(z) = \sum_{n=1}^\infty \widetilde{b}_n(t)z^n.  
\end{split}
\end{equation}

The transformations (\ref{eq611}) do not preserve the time-one probability measure $\mu_1$ in generic cases. For the reader's convenience, 
we state the results in the two cases of monotone and c-monotone convolutions separately. 
\begin{thm} \label{thm112} (Monotone case)  
(1) If $\nu_1 = \widetilde{\nu}_1$ is not a delta measure, then there exists an integer $n$ such that 
\begin{equation}\label{eq212}
\widetilde{B}_2 = \Big(1+ \frac{2 \pi i n}{s_1}\Big)B_2. 
\end{equation}

(2) If $\nu_1= \widetilde{\nu}_1$ is a delta measure, then there exists an integer $n$ such that 
\begin{equation}\label{eq213}
\widetilde{B}_2 = B_2 + 2 \pi i n. 
\end{equation}   
\end{thm}

\begin{thm}\label{thm113} (C-monotone case) 
(1) If $\nu_1= \widetilde{\nu}_1$ is not a delta measure, then $B_2$ and $\widetilde{B}_2$ satisfy the relation (\ref{eq212}). 
Moreover, there exists an integer $m$ such that 
\begin{equation}\label{eq663}
\widetilde{B}_1 = 2 \pi i m - \frac{2 \pi i n r_1}{s_1} + \Big(1 + \frac{2 \pi i n}{s_1}\Big)B_1. 
\end{equation}

(2) If  $\nu_1 =  \widetilde{\nu}_1$ is a delta measure, then $B_2$ and $\widetilde{B}_2$ satisfy the relation (\ref{eq213}) for an $n$. 
$B_1$ and $\widetilde{B}_1$ necessarily satisfy $\widetilde{r}_1 \in r_1 + 2\pi i \mathbb{Z}$. In addition, there are several cases.  

(a) If $s_1 \in i \real \cap (2 \pi i \rat)^c$, then there exists an integer $m$ such that 
(\ref{eq663}) holds. 

(b) We assume that $s_1 = 2 \pi i \frac{p}{q}$ for  integers $p \neq 0$ and $q >0$ which cannot be divided by a common prime number. We moreover assume that 
$\widetilde{s}_1 \neq 0$ if $s_1 \in 2\pi i \mathbb{Z}$. Then $\widetilde{r}_{j+1} = r_{j+1}$ for $j \in \nat \cap (q\nat)^c$, where $\nat = \{1,2,\cdots \}$. There are no restrictions on $r_j$ and $\widetilde{r}_j$
 for $j \in q\nat +1$.  

If one of $s_1$ and $\widetilde{s}_1$ is $0$, there are three cases. 

 (c) If $s_1 = 2 \pi i p$ for a non-zero integer $p$ and $\widetilde{s}_1 = 0$, then  
 $\widetilde{r}_j = 0$ for $j \geq 2$. There are no restrictions on $r_j$ for $j \geq 2$. 

(d) If $s_1 = 0$ and $\widetilde{s}_1 =  2 \pi i q$ for a non-zero integer $q$, then  
 $r_j = 0$ for $j \geq 2$. There are no restrictions on $\widetilde{r}_j$ for $j \geq 2$. 
 
(e) If $s_1 = \widetilde{s}_1 =0$, then $r_j = \widetilde{r}_j $ for $j \geq 2$.  
%
\end{thm}

It is difficult to formulate the above theorems in terms of transformations. Let us focus on the monotone case. 
In the case (1), the transformations 
\begin{equation}
 B_2 \mapsto   \Big(1 + \frac{2 \pi i n}{s_1}\Big)B_2  
\end{equation}
preserve the time one probability measure. However,  $\Big(1 + \frac{2 \pi i n}{s_1}\Big)B_2$ may not map the unit disc $\disc$ into 
the left half plane $\{z \in \comp: \re z \leq 0 \}$. We take the function $B_2(z) = z-a$ for $\re a \geq 1$ as an example. 

If $a =1$, then $B_2(z)$ is equal to $z-1$ whose image is tangent to the imaginary axis in the complex plane.  If this image is rotated by however small angle, it has a nonempty 
intersection with the imaginary axis. Therefore, 
the image of $(1 - 2 \pi i n)B_2$ never be contained in the left half plane. This implies that 
there is no other function $B_2$ which generates the same measure at time one.  

Next let $a$ be a sufficiently large real number. Then we can easily prove that the function  $\Big(1 - \frac{2 \pi i n}{a}\Big)B_2$ maps 
$\disc$ into the left half plane for some non-zero integer $n$. This means that the function $B_2$ is not unique for the time one measure $\mu_1$.

To prove the two theorems, we need to express $s_n$ and $r_n$ in terms of $b_n(1)$ and $a_n(1)$. 
This is done through the following lemmata. 

\begin{lem}\label{lem2211} Let  $f_n(t):= a_n(t)e^{-r_1 t}$ and $g_n(t):= b_n(t)e^{-s_1t}$ for $n \geq 1$.  
There exist polynomials $P_{(l_1,\cdots,l_n; k_1,\cdots,k_n)}(x)$ and $Q_{(k_1,\cdots,k_n)}(x)$  for $n \geq 1$ and $k_i,l_i \geq 0$ such that 
\begin{align}
&f_n(t) = \frac{r_n}{r_1}\frac{e^{(n-1)s_1 t}-1}{n-1} + \sum_{\substack{1 \leq l_1 + \cdots + l_{n-2} + k_1+ \cdots + k_{n-2}  \leq n-1,\\ l_i \geq 0,~ k_i \geq 0}} \prod_{j=2}^{n-1}\Big( \frac{r_j}{s_1}\Big)^{l_j}\prod_{j=2}^{n-1} \Big( \frac{s_j}{s_1}\Big)^{k_j} P_{(l_1,\cdots,l_n; k_1,\cdots,k_n)}(e^{s_1t}),\\
& g_n(t) = \frac{s_n}{s_1} \frac{e^{(n-1)s_1 t}-1}{n-1} +  \sum_{\substack{1 \leq k_1+ \cdots + k_{n-2} \leq n-1,\\ k_i \geq 0}}  \prod_{j=2}^{n-1}\Big( \frac{s_j}{s_1}
\Big)^{k_j} Q_{(k_1,\cdots,k_{n-2})}(e^{s_1t})
\end{align}
for $n \geq 2$. The summations are understood to be $0$ for $n=2$. If $n=1$, $f_1(t) = g_1(t) = 1$. 
\end{lem}
\begin{proof}
From the coefficients of $z^n$ in the differential equations (\ref{eq51}) and (\ref{eq52}), it holds that   
\begin{align}
&\frac{d}{dt}a_n (t) = r_1a_n(t) + r_n a_1(t) b_1(t)^{n-1} + \sum_{m=2}^{n-1}r_m \sum_{\substack{ k_1+ \cdots + k_{m} = n,\\ k_i \geq 1}} b_{k_1}(t)\cdots b_{k_{m-1}}(t)a_{k_{m}}(t), \\
&\frac{d}{dt}b_n (t) = s_1b_n(t) + s_n b_1(t)^{n} + \sum_{m=2}^{n-1}s_m \sum_{\substack{k_1+ \cdots + k_{m} = n,\\ k_i \geq 1}} b_{k_1}(t)\cdots b_{k_{m}}(t)
\end{align}
for $n \geq 2$. If $n = 2$, the summations are understood to be $0$. $\frac{d}{dt}a_1(t) = r_1a_1(t)$ and $\frac{d}{dt}b_1(t) =s_1b_1(t)$ for $n=1$.  We note that initial conditions are $a_1(0) = b_1(0) =1$, $a_n(0)=b_n(0) = 0$ for 
$n \geq 2$. 
In terms of $f_n(t)$ and $g_n(t)$, we have 
\begin{align}
&\frac{d}{dt}f_n (t) =  r_n  b_1(t)^{n-1} + \sum_{m=2}^{n-1}r_m \sum_{\substack{k_1+ \cdots + k_{m} = n,\\ k_i \geq 1}} b_{k_1}(t)\cdots b_{k_{m-1}}(t)f_{k_{m}}(t), \\
&\frac{d}{dt}g_n (t) = s_n b_1(t)^{n-1} + \sum_{m=2}^{n-1}s_m \sum_{\substack{k_1+ \cdots + k_{m} = n,\\ k_i \geq 1}} b_{k_1}(t)\cdots b_{k_{m-1}}(t)g_{k_{m}}(t). 
\end{align}
Then the claim can be proved by induction. 
\end{proof}

\begin{lem}\label{lem21}Let $a_n:=a_n(1)$ and $b_n:=b_n(1)$. Then $a_1 = e^{r_1}$ and $b_1= e^{s_1}$. 
We assume that $s_1 \notin 2\pi i\rat$. 
Then for $n \geq 2$, there exist polynomials $P_n^{(1)}(x_1,\cdots,x_{n-1},y_1,\cdots,y_{n-1})$, $Q_n^{(1)}(x, y)$, $P_n^{(2)}(x_1,\cdots,x_{n-1})$, $Q_n^{(2)}(x)$ such that 
\begin{align}
&\frac{r_n}{s_1} = \frac{n-1}{a_1(b_1^{n-1}-1)}a_n + \frac{P^{(1)}_n(a_1,\cdots,a_{n-1},b_1,\cdots,b_{n-1})}{Q_n^{(1)}(a_1,b_1)}, \\
&\frac{s_n}{s_1} = \frac{n-1}{b_1(b_1^{n-1}-1)}b_n + \frac{P^{(2)}_n(b_1,\cdots,b_{n-1})}{Q_n^{(2)}(b_1)}. 
\end{align}
$Q^{(1)}_n$ and $Q^{(2)}_n$ satisfy that $Q_n^{(1)}(x,y) \neq 0$ for $x \neq 0, y \notin \tor$ and $Q_n^{(2)}(x) \neq 0$ for $x \notin \tor \cup \{0\}$.  
Therefore $Q_n^{(1)}(a_1,b_1) \neq 0$ and $Q_n^{(2)}(b_1) \neq 0$ under the assumption  $s_1 \notin 2\pi i\rat$. 
\end{lem}
\begin{proof}
This claim can be proved by a simple argument of induction and by Lemma \ref{lem2211}. 
\end{proof}
\begin{proof}[Proof of  the theorems] 
Since $a_1 = \widetilde{a}_1$ and $b_1 = \widetilde{b}_1$, immediately $\widetilde{r}_1 \in r_1 + 2\pi i \mathbb{Z}$ and $\widetilde{s}_1 \in s_1 + 2\pi i \mathbb{Z}$ follow. 
Therefore if $s_1 \in i \real$, $\widetilde{s}_1$ also belong to $i\real$. The proof of Proposition \ref{prop21} claims that $s_k = \widetilde{s}_k =0$ for $k \geq 2$. Thus we have proved Theorem \ref{thm112} (2). 
We next assume that $s_1 \notin 2\pi i \rat$. $\mu_1 = \widetilde{\mu}_1$ and $\nu_1 = \widetilde{\nu}_1$ are equivalent to $a_n = \widetilde{a}_n (:=\widetilde{a}_n(1))$ and 
$b_n = \widetilde{b}_n (:=\widetilde{b}_n(1))$. By Lemma \ref{lem21}, these conditions are also equivalent to  $\frac{s_k}{s_1}=\frac{\widetilde{s_k}}{\widetilde{s_1}}$ and $\frac{r_k}{s_1}=\frac{\widetilde{r_k}}{\widetilde{s_1}}$ 
for $k \geq 2$. Therefore, there exist integers $m,n$ such that $\widetilde{s_k}= \Big( 1+\frac{2\pi i n}{s_1}\Big) s_k$ for $k \geq 1$, $\widetilde{r}_1 = r_1 + 2\pi i m$ and
$\widetilde{r_k}= \Big( 1+\frac{2\pi i m}{s_1}\Big) r_k$ for $k \geq 2$. Thus we have proved Theorem \ref{thm112} (1) and Theorem \ref{thm113} (1), (2-a). 
Theorem \ref{thm113} (2-b)-(2-e) can be proved by applying Proposition \ref{prop31}.
\end{proof}


\subsection{Connections to the Boolean convolution}
We would like to construct $\dismrhd_c$-convolution semigroups similarly to the additive case. First we 
define a multiplicative version of the $t$-transformation introduced in \cite{BW1}.  
However, multiplicative Boolean infinite divisibility does not hold for all probability measures. Fortunately, 
we have the following Theorem \ref{thm11}, so that the multiplicative $t$-transformation of a $\dismrhd$-convolution semigroup 
can be defined. 
We denote by $ID(\dismrhd, \tor)$ the set of infinitely divisible distributions on $\tor$ for the multiplicative monotone convolution. 
\begin{thm}\label{thm11} Let $(\mu, \nu) \in \pa(\tor) \times \pa(\tor)$ be $\dismrhd_c$-infinitely divisible. Then both $\mu$ and $\nu$ belong to $ID(\putimes, \tor)$. In particular, we have $ID(\dismrhd, \tor) \subset ID(\putimes, \tor)$. 
\end{thm}
\begin{proof} 
We may assume that $\mu, \nu \neq \omega$; otherwise, the claim is trivial. Let $\{(\mu_t, \nu_t) \}_{t \geq 0}$ be a weakly continuous $\dismrhd_c$-convolution semigroup with $(\mu_0, \nu_0) = (\delta_1, \delta_1)$ and $(\mu_1, \nu_1) = (\mu, \nu)$.  
By Theorem \ref{thm4}, there exists a weakly continuous $\dismrhd_c$-convolution semigroup $\{(\mu_t, \nu_t) \}_{t \geq 0}$ with $(\mu_0, \nu_0) = (\delta_1, \delta_1)$ and $(\mu_1, \nu_1) = (\mu, \nu)$.  Then Theorem \ref{thm2121} enables us to take two vector fields $B_1$ and $B_2$ defined in $\disc$ such that 
\[
\begin{split}
&\eta_{\mu_t}(z)=z\exp\Big(\int_0 ^t B_1(\eta_{\nu_s}(z))ds\Big),  \\
&\eta_{\nu_t}(z) = z\exp\Big(\int_0 ^t B_2(\eta_{\nu_s}(z))ds\Big). 
\end{split}
\]
These expressions imply that $\frac{\eta_{\mu}(z)}{z}$ and  $\frac{\eta_{\nu}(z)}{z}$ do not have zero points. 
\end{proof}
\begin{rem}
It may be a nontrivial question whether the same relation holds in the c-free case. 
\end{rem}

When $\mu \in ID(\putimes, \tor)$, $\frac{\eta_\mu(z)}{z}$ is written as $e^{u_\mu(z)}$. The representation is unique if we impose the condition $\im u_\mu(0) \in [0, 2\pi)$, for instance (see \cite{Fra2}). We always choose this branch and then we define $\mu^{\hutimes t}$ by $\eta_{\mu^{\hutimes t}}(z) = ze^{tu_\mu(z)}$.
If we try to define Boolean convolution semigroups, this ambiguity of the branches necessarily occurs. As a result,   
$\dismrhd_c$-convolution semigroups cannot be simply constructed  from $\putimes$-convolution semigroups, as we see below.   We remark that the relation 
$\mu^{\hutimes s} \putimes \mu^{\hutimes t} = \mu^{\hutimes s + t}$ holds for all $s, t \geq 0$, but $\big{(}\mu^{\hutimes s}\big) ^{\hutimes t}$ is not equal to $\mu^{\hutimes st}$ for general $s, t \geq 0$. 

\begin{defi}
(1) We define a map $\mathcal{V}_t: ID(\putimes, \tor) \to ID(\putimes, \tor)$ by $\mathcal{V}_t(\mu):= \mu^{\hutimes t}$ for $t \geq 0$. \\
(2) We define a map $\Theta^{u, v}: ID(\putimes, \tor) \times ID(\putimes, \tor) \to ID(\autimes, \tor)$ by 
$\Theta^{u, v}(\mu, \nu):= \mu^{\hutimes u} \putimes \nu^{\hutimes v}$ for $u, v \geq 0$. 
\end{defi}


\begin{defi}\label{defi3}
(1) Let $\{ \nu_t \}_{t \geq 0}$ be a weakly continuous $\dismrhd$-convolution semigroup generated from a vector field  $B_2 ^{\nu}$.  
We define $\{ (\mu^r _t, \nu_t) \}_{t \geq 0}$ for $r \geq 0$ by the pair of vector fields $(rB_2 ^{\nu}, B_2 ^\nu)$. 

(2) Let $\{ (\kappa_t, \lambda_t) \}_{t \geq 0}$ and $\{ (\nu_t, \lambda_t) \}_{t \geq 0}$ be weakly continuous $\dismrhd_c$-convolution semigroups 
generated respectively from $(B_1 ^{\kappa, \lambda},B_2^\lambda)$ and $(B_1^{\nu,\lambda},B_2^\lambda)$. 
We define $\{ (\mu^{u,v}_t, \lambda_t) \}_{t \geq 0}$ by the pair of vector fields $(uB_1^{\kappa, \lambda} + v B_1 ^{\nu, \lambda}, B_2 ^\lambda)$. 
These definitions are parallel to the additive c-monotone convolution semigroups \cite{Has3}. 
\end{defi}
The definitions of $\mu^r _t$ and $\mu^{u,v}_t$ are identical to $\mathcal{V}_r(\mu_t)$ and $\Theta^{u,v}(\kappa_t, \nu_t)$, respectively, for $r \in \nat$. For general $r > 0$, however, they are not. We have the following properties for small $r > 0$.  
\begin{prop} In the above notation, we have the following. \\
(1) If $t$ satisfies $0 \leq t\,\im B_2 ^\nu(0) < 2\pi $, we have $\mathcal{V}_r (\nu_t) = \mu^r _t$. \\
(2) If $u, v$ satisfy $0 \leq u\,\im B_1^{\kappa, \lambda}(0) < 2\pi $ and $0 \leq v \,\im B_1 ^{\nu, \lambda}(0) < 2\pi$, 
we have $\mu^{u,v} _t = \Theta^{u,v}(\kappa_t, \nu_t)$. 
\end{prop}
\begin{proof}
We only prove (1) since (2) can be proved in the same way. By Theorem \ref{thm11} we can write $\eta_{\nu_t}(z) = z e^{u_t(z)}$ for an analytic function $u_t$ satisfying $u \in C^\omega([0, \infty) \times \disc)$ and $ u_0(z) = 0$. The differential equation for $\eta_{\nu_t}$ becomes $\frac{d}{dt}u_t(z) = B_2 ^\nu(ze^{u_t(z)})$. Then we obtain 
$\im u_t(0) = t \,\im B_2 ^\nu (0)$. By definition, $\eta_{\mu^r _t}(z) = ze^{ru_t(z)}$ and $\eta_{\mathcal{V}_r (\nu_t)}(z) = ze^{ru_t(z)}$, the latter of which holds when $0 \leq t\,\im B_2 ^\nu(0) < 2\pi $. 
\end{proof}

In view of these results, the multiplicative version of $t$-transformation, which means the time evolution with respect to the Boolean convolution, does not work  so 
 well in comparison with the additive convolution. This problem needs to be investigated further, including the case of free and c-free convolutions.

We show examples where explicit forms of probability measures can be calculated.   
\begin{exa}\label{exa3}
The relation $\frac{1 + \eta_\mu(z)}{1 - \eta_\mu (z)} = \int_{-\pi} ^{\pi}\frac{e^{i\theta} +z}{e^{i\theta}-z}\mu(-d\theta)$ is useful in the following calculations. 
Moreover, Theorem \ref{thm6} is also convenient. \\
(1) If $B_2(z) = -a +ib = ib - a\int_{-\pi} ^\pi \frac{e^{i\theta} +z}{e^{i\theta}-z}\omega(d\theta)$ ($a > 0, b \in \real$), we have $\eta_{\nu_t}(z) = ze^{(-a +bi)t}$ and 
\[
\nu_t(d\theta) = \frac{1}{2\pi}\frac{1 - e^{-2at}}{1 + e^{-2at}-2e^{-at}\cos(\theta - bt)}d\theta. 
\]
This is identical to the density of the Poisson kernel. 
$\mu^r _t$ in Definition \ref{defi3} (1) is obtained only by the transformations $a \mapsto ra$ and $b \mapsto rb$ since $B_2$ is 
a constant. \\
(2) If $B_2(z) = a(z-1)= -a\int_{-\pi} ^\pi \frac{e^{i\theta} +z}{e^{i\theta}-z}(1-\cos \theta )\omega(d\theta)$ ($a > 0$), then we have $\eta_{\nu_t}(z) = \frac{z}{(1-e^{at})z + e^{at}}$ and 
\[
\nu_t = (1 - e^{-at})\omega + e^{-at}\delta_1. 
\] 
We can easily check that $\frac{\eta_{\nu_t}(z)}{z}$ does not have a zero point.  
 $\mu^r _t$ is obtained by the equation $\frac{d}{dt}\log\eta_{\mu^r _t}(z) = ra(\eta_{\nu_t}(z) - 1)$ and we have 
 \[
\eta_{\mu^r _t}(z) = z (z + (1 - z)e^{at})^{-r}. 
\]
\end{exa}

\section*{Acknowledgement} 
The author thanks Professor Izumi Ojima for discussions on multiplicative convolutions and complex analysis. 
He also thanks Professor Uwe Franz for suggesting the contents of Subsection \ref{subsec5}. 
This work was supported by JSPS (KAKENHI 21-5106) and Global COE Program at Kyoto University.

\providecommand{\bysame}{\leavevmode\hbox to3em{\hrulefill}\thinspace}

\end{document}